\documentclass[10pt]{amsart}
\usepackage{amsmath,amssymb,enumerate}
\usepackage{hyperref}
\newtheorem{thm}{Theorem}[section]
 \numberwithin{equation}{section} 
 \numberwithin{figure}{section} 
 \theoremstyle{plain}
 \theoremstyle{plain}    
 \newtheorem{cor}[thm]{Corollary} 
 \theoremstyle{plain}    
 \newtheorem{prop}[thm]{Proposition} 
 \newtheorem{defi}[thm]{Definition}
 \theoremstyle{plain}    
 \newtheorem{lem}[thm]{Lemma} 
 \theoremstyle{remark}
 \newtheorem{rem}[thm]{Remark}
 \theoremstyle{definition}

\begin{document}

\title[Stability for MA equations]{$L^1$-Stability for complex Monge--Ampère equations}
\author{Songchen Liu}
\author{Liyou Zhang}
\address{School of mathematical sciences, Capital Normal University, Beijing, China}
\email{{2230502073@cnu.edu.cn}}
\email{zhangly@cnu.edu.cn}
\date{\today}
\begin{abstract}
    We first establish the weak stability results for solutions of complex Monge--Ampère equations in relative full mass classes, extending the results known to hold in the full mass class. 
    Building on weak stability, we then prove the $\mathcal{C}^{k,\alpha}$ stability of solutions to complex Monge–Ampère equations on quasi-projective varieties. 
    As an application, we study the limit of the singular Ricci-flat metrics on $\mathbb{Q}$-Calabi--Yau projective varieties, inspired by Tosatti's work on Calabi–Yau projective manifolds.
\end{abstract}

\maketitle

\section{Introduction}

Let $(X,\omega)$ be a compact Kähler manifold of complex dimension $n$, and let $\theta$ be a closed real $(1,1)$-form. 
The systematic study of complex Monge–Ampère equations in big cohomology classes has been initiated in \cite{BEGZ10}. When $\{\theta\}$ is a big cohomology class (i.e. it contains a Kähler current), it has been shown that there exists a unique solution $\varphi\in {\rm PSH}(X,\theta)$ with full Monge--Ampère mass (i.e. $\int_X \theta_\varphi^n = \int_X \theta_{V_\theta}^n$) such that
\begin{equation}\label{eq 1}
    \theta_\varphi^n = \mu,~~\sup_X \varphi = 0,
\end{equation}
if and only if $\mu$ is a normalized non-pluripolar measure on X, where $\theta_\varphi^n := \left( \theta+ dd^c \varphi \right)^n$ is non-pluripolar product. 

  Recently, the stability of solution of (\ref{eq 1}) has attracted lots of attentions. In \cite[Proposition A]{1}, Guedj and Zeriahi proved the following interesting results, which generalized Cegrell--Kolodziej's result \cite{CK06} in the local setting.\par
  Let $\{ \theta \} $ be a fixed big cohomology class, and let $\mu_j = \theta_{\varphi_j}^n,j=1,....$ be non-pluripolar measures, where $\varphi,\varphi_j \in \mathcal{E}(X,\theta)$ (i.e. $\int_X \theta_{\varphi_j}^n =\int_X \theta_\varphi^n = \int_X \theta_{V_\theta}^n$) and $\sup_X \varphi = \sup_X \varphi_j = 0$. If $\|\mu_j - \mu \| \to 0$, then 
  $$
  \| \varphi_j -\varphi \|_{L^1(\omega^n)}\to 0.
  $$
  Here $\|\mu_j - \mu\|$ denotes the total variation of the signed measure $\mu_j - \mu$.\par
  
 In \cite[Theorem 5.17]{DDL23}, Darvas, Di Nezza, and Lu proved that when $\{ \theta \}$  is a big cohomology class, for each $\theta$-model potential $\phi$ (see Definition \ref{def 1}) and normalized non-pluripolar measure $\mu$, there exists a unique $\varphi \in \mathcal{E}(X, \theta, \phi)$ (see Definition \ref{def 2}) and $\sup_X \varphi =0$ such that $ \varphi$ is the solution of (\ref{eq 1}).\par
 
 In this paper, we investigate the stability of solution of (\ref{eq 1}) in the relative full mass class $\mathcal{E}(X,\theta,\phi)$. 
 More precisely, assume that that $\{\theta\},\{\theta^j\}$ for $j=1,2,...$, are big cohomology classes, $\phi_j~({\rm resp.}~\phi)$ is a $\theta^j$-model potential (resp. $\theta$-model potential), and $\mu_j,\mu$ are normalized non-pluripolar measures such that
 \begin{equation}\label{eq 2}
 \begin{split}
    (\theta^j_{\varphi_j})^n = \mu_j, ~\varphi_j \in \mathcal{E}(X,\theta^j,\phi_j)&~{\rm and}~\sup_X \varphi_j  =0; j=1,2,...,\\
 \theta_\varphi^n = \mu,  ~\varphi \in \mathcal{E}(X,\theta,\phi)&~{\rm and}~\sup_X \varphi  =0.
 \end{split}
 \end{equation}
 Our goal is to find some \emph{suitable} conditions on  $\theta^j,\phi^j,\mu_j$ such that $\varphi_j \to_{L^1(\omega^n)} \varphi$.\par
  
   We denote $\theta^j \to_{\mathcal{C}^0} \theta$ (resp. $\theta^j \to_{\mathcal{C}^+} \theta$) if for all $\epsilon>0$, there exists $j_0>0$ such that 
  $$
  -\epsilon \omega \leq \theta_j - \theta \leq \epsilon \omega~({\rm resp.}~ \theta_j - \theta \leq \epsilon \omega)
  $$
  for all $j>j_0$. Clearly, $\mathcal{C}^0$-convergence implies $\mathcal{C}^+$-convergence. We establish the following $L^1$-stability for (\ref{eq 2}).
 \begin{thm}\label{Proposition A}$($=Theorem \ref{Proposition A'}$)$
   Assume that $\theta^j \to_{\mathcal{C}^+} \theta$, and $\phi_j \to_{L^1(\omega^n)} \phi$. If $\| \mu_j - \mu \| \to 0$, then
$$
\| \varphi_j -\varphi \|_{L^1(\omega^n)} \to 0.
$$
\end{thm}

\begin{rem}
    In particular, if we take $\theta^j = \theta$ and  $\phi_j = \phi = V_{\theta}$,
    this special case of Theorem \ref{Proposition A} recovers Guedj--Zeriahi's result, \cite[Proposition A]{1}.
\end{rem}

If $\phi$ dominates all $\phi_j$ for $j=1,2,...$, the condition of $\phi_j \to_{L^1(\omega^n)} \phi$ is superfluous.

\begin{cor}\label{Proposition B}$($=Corollary \ref{Proposition B'}$)$
   Assume that $\theta^j \to_{\mathcal{C}^+} \theta$, and $\phi_j \leq \phi$. If $\| \mu_j - \mu \| \to 0$, then
$$
\| \varphi_j -\varphi \|_{L^1(\omega^n)} \to 0.
$$
\end{cor}

If we set $\phi = V_\theta$, then for arbitrary $\theta^j$-model potentials $\phi_j $, we have
\begin{cor}\label{Proposition C}$($=Corollary \ref{Proposition C'}$)$
   Assume that $\theta^j \to_{\mathcal{C}^+}\theta$, and $\phi = V_\theta$. If $\| \mu_j - \mu \| \to 0$, then
$$
\| \varphi_j -\varphi \|_{L^1(\omega^n)} \to 0.
$$
\end{cor}

    The main ingredient in the proof of the above $L^1$-stability results is as follows. 
    The stability of subsolutions of complex Monge–Ampère equations, established by Darvas, Di Nezza, and Lu \cite[Lemma 5.16]{DDL23}, asserts that subsolutions of these equations are stable under $L^1$-convergence and has been further generalized in Lemma \ref{lem 3.2}.
    We then construct a new family of complex Monge--Ampère equations, for which $\varphi_j$ are subsolutions. 
    Applying the weak compactness of quasi-psh functions \cite[Proposition 8.5]{GZ17}, and the uniqueness principle \cite[Theorem 3.13]{DDL23}, we have the $L^1$-limit of $\varphi_j$ is identical to $\varphi$.\par

    There have also been some interesting results regarding the stability of (\ref{eq 2}).
    In \cite[Section 3]{DDL21}, Darvas, Di Nezza, and Lu introduced the distance of singularity types, $d_{\mathcal{S}}$, which plays a crucial role in the stability of (\ref{eq 2}).
    In \cite[Theorem 1.4]{DDL21}, they dealt with the case that $\theta^j = \theta$, $\mu_j,\mu$ has $L^p(\omega^n)$ density, and $\phi_j \to_{d_{\mathcal{S}}} \phi$, concluding that $\varphi_j \to \varphi$ in capacity. 
    In \cite{DV22,DV24}, Do and Vu considered the stability of (\ref{eq 2}) under the conditions that $\theta^j \to_{\mathcal{C}^0} \theta$, $\phi_j \to_{d_{\mathcal{S},(A+1)\omega}} \phi$ (where $A\omega \geq \theta^j,\theta$), and $\| \mu_j - \mu \| \to0$, and show that $\varphi_j \to_{d_{\rm cap}} \varphi$. Moreover, Do--Vu even  established a quantitative stability result, see \cite[Theorem 1.4]{DV22}.\par
    
    In both cases, Darvas--Di Nezza--Lu and Do--Vu established that $\varphi_j \to \varphi$ in a norm induced by capacity, which is stronger than $L^1$-convergence. 
    Note that the $d_\mathcal{S}$-convergence of model potentials implies $L^1$-convergence when cohomology classes are fixed, as shown in \cite[Theorem 5.6]{DDL21}.  However, the converse does not generally hold. 
    Therefore, we generalize our approach in Theorem \ref{Proposition A} by adopting $L^1$-convergence of model potentials instead. 
    As a result, this ensures only $L^1$-convergence of solutions, which nevertheless suffices for our subsequent results.\par
    
    As an application of Theorem \ref{Proposition A}, we have the $L^1$-compactness of solutions to (\ref{eq 2}). 
    Then, by Zeriahi's uniform Skoda integrability theorem \cite{Ze01}, we can study the $\mathcal{C}^{k,\alpha}$-convergence of solutions to complex Monge--Ampère equations on quasi-projective varieties, based on a result due to Di Nezza and Lu, \cite{DNL14}.\par

  Let $\theta$ be a smooth, closed semi-positive real $(1,1)$-form on $X$ such that $\int_X \theta^n >0$. 
  Assume $\psi^{\pm} \in {\rm QPSH}(X) \cap \mathcal{C}^\infty(X\backslash D)$, where $D$ is a divisor on $X$, and $e^{- \psi^-} \in L^1(\omega^n)$. 
  Consider the following complex Monge--Ampère equations:
 \begin{equation}\label{eq 1.4}
 \theta_\varphi^n = ce^{\psi^+ - \psi^-} \omega^n,~\varphi \in \mathcal{E}(X,\theta),
 \end{equation}
 where $c$ is a normalized constant. By Guan–Zhou's theorem \cite{GZh15}, it follows that $e^{-\psi^-} \in L^p(\omega^n)$ for some $p > 1$. Consequently, using Guedj--Zeriahi's $L^\infty$ estimates for semi-positive form \cite[Theorem 1.3]{GZ07} (see also \cite{BEGZ10, DDL18, GL21} for more general cases), we obtain that $\varphi \in L^\infty(X)$.\par

  The equation (\ref{eq 1.4}) is known as {\it the complex Monge--Ampère equation on quasi-projective varieties}, and it has been extensively studied during the past years.
  For instance, Eyssidieux--Guedj--Zeriahi \cite{EGZ09}, Pǎun \cite{Pa08} studied the regularity of solutions on $X \backslash D$ for such equations when $\psi^\pm$ have analytic singularities, and investigated the singular Kähler--Einstein metrics on normal varieties.

  In \cite{DNL14}, Di Nezza and Lu proved that the solution $\varphi$ of (\ref{eq 1.4}) is smooth outside $D\cup E$, where $E$ is an effective {simple normal crossing} (s.n.c.) $\mathbb{R}$-divisor on $X$ such that $\{ \theta \} - c_1(E)$ is a Kähler class. 
   Now, applying Theorem \ref{Proposition A}, we present a $\mathcal{C}^{k,\alpha}$ stability result for (\ref{eq 1.4}):
\begin{thm}\label{thm D}$($=Theorem \ref{thm D'}$)$
  Assume $\theta^j $ are smooth, closed semi-positive real $(1,1)$-forms on $X$ such that $\theta^j \to_{\mathcal{C}^0} \theta$. Then, for $j>>1$, the normalize solution $  \varphi_j $ of the complex Monge--Ampère equation
  \begin{equation}\label{eq 1.3}
       (\theta^j_{\varphi_j})^n = c_j e^{\psi^+ - \psi^-} \omega^n,~\varphi_j \in \mathcal{E}(X,\theta^j)
  \end{equation}
    is smooth outside $D\cup E$ for each $j$. If the local potentials of $\theta^j$ are $\mathcal{C}^{k,1}(U)$-convergent to the local potential of $\theta$ for any coordinate neighborhood $U \Subset X \backslash (D \cup E)$ and $\forall k \in \mathbb{Z}^+$, then there exists $0<\alpha<1$ such that
    $$
    {\varphi_j} \to_{\mathcal{C}^{k,\alpha}_{loc}(X \backslash (D\cup E))} \varphi,~\forall  k \in \mathbb{Z}^+.
    $$
\end{thm}

Actually, we can consider $\psi_j^\pm$ as a sequence of quasi-psh functions, rather than fixing $\psi^\pm$ as in Theorem \ref{thm D}, while assuming that $e^{\psi_j^+ - \psi_j^-} \to_{L^1(\omega^n)} e^{\psi^+ - \psi^-}$. Under this assumption, the same conclusion in Theorem \ref{thm D} holds; see Section~\ref{sec 4.3}.
 
Let $D$ be an effective s.n.c. $\mathbb{R}$-divisor on $X$.
Under certain positivity assumptions on $K_X + D$, related works on stability of complex Monge–Ampère equations on quasi-projective varieties, have also been studied recently, c.f. \cite{BG22,DV23}.
Among these, Biquard, Guenancia primarily focused on $\mathcal{C}^\infty_{loc}$-convergence \cite{BG22}; Dang and Vu, on the other hand, investigated weak convergence \cite{DV23}.
Their geometric motivation comes from the study of {\it Kähler metrics with singularities along a divisor}.\par


\subsection{Organization} In Section \ref{sec 2}, we briefly recall some material from pluripotential theory, which is developed in full detail in \cite{Bou04, BEGZ10, DDL23,GZ17}. 
In Section \ref{sec 3}, we establish several weak stability results for (\ref{eq 2}). 
In Section \ref{sec 4}, we first revisit key results from \cite{DNL14}, and then prove Theorem \ref{thm D}. 
Finally, in Section \ref{sec 4.3}, we provide a slight generalization of Theorem \ref{thm D} to the case of variations in the right-hand side of (\ref{eq 1.3}); and in section \ref{sec 4.4}, we provide an application of our main results  to Calabi–Yau varieties.

\subsection{Acknowledgment} 
We are grateful to T. Darvas for his insightful comments, and for pointing out an error in Lemma \ref{lem 3.2}, and suggesting a concise approach, which greatly enhanced the article.
We also thank Shiyu Zhang for fruitful and helpful discussions about section \ref{sec 4.4}.
The first author is indebted to T. Darvas for his invaluable and stimulating discussions over an extended period.
 This work is partially supported by NSFC grants 12071310 and 12471079.

\section{Preliminaries}\label{sec 2}

  Let $(X, \omega)$ be a compact Kähler manifold of complex dimension $n$. Let $\theta$ be a closed real $(1,1)$-form on $X$ and $\varphi$ be a upper semi-continuous (u.s.c.) function, $\varphi: X \to \mathbb{R}\cup\{ -\infty \}$. We denote $\varphi \in {\rm PSH}(X,\theta)$ {\it iff} $\varphi \in L^1(\omega^n)$ and
$$
\theta_\varphi:=\theta+dd^c \varphi\geq 0
$$
 in the sense of currents, where $dd^c:=\frac{\sqrt{-1}}{\pi}\partial\bar{\partial} $.
 A u.s.c.~function $\varphi$ is said to be a quasi-psh function if it can be written as the sum of a psh function and a smooth function locally. We denote by $ \mathrm{QPSH}(X)$ the set of quasi-psh functions on $X$. Clearly, we have ${\rm PSH}(X, \theta) \subset \mathrm{QPSH}(X)$.\par
 The cohomology class of a  smooth closed real (1,1)-form $\theta$ is said to be {\it big} if there exists $\varphi\in {\rm PSH}(X,\theta)$ such that $\theta_\varphi$ dominates a Kähler form. Unless otherwise stated, the cohomology classes of $\theta$ and $\theta^i$ are all big throughout this section.\par


 Let $\varphi,\psi \in {\rm PSH}(X,\theta)$. Then $\psi$ is said to be {\it less singular} than $\varphi$, say $\varphi \preceq  \psi$, 
 if they satisfy $\varphi \leq \psi + C$ for some $C \in \mathbb{R}$. We say that $\varphi$ has the {\it same singularity} as $\psi$, say $\varphi\backsimeq  \psi$,
 if $\varphi\preceq \psi$ and $\psi \preceq\varphi$. Let $f$ be a real-valued function on $X$. We define the {\it envelope} of $f$ in the class ${\rm PSH}(X, \theta)$ by
$$
 P_\theta(f) := ( \{\sup{ u \in {\rm PSH}(X,\theta) : u \leq f} \})^*.
$$
  Set $V_\theta: = P_\theta(0)$. Then $V_\theta$ has {\it minimal singularity}, i.e. $V_\theta \succeq \varphi$, $\forall \varphi \in {\rm PSH}(X,\theta)$.\par

  \subsubsection{Non-pluripolar product} The ample locus ${\rm Amp}(\{\theta\})$ of $\theta$ is the set of points $x \in X$ such that there exists a Kähler current $T \in \{\theta\}$ with analytic singularity type and smooth in a neighborhood of $x$. The ample locus ${\rm Amp}(\{\theta\})$ is a  nonempty Zariski open subset, see \cite{Bou04}. Obviously, one has $V_\theta \in L^\infty_{loc}({\rm Amp}(\{ \theta \}) )$.\par
  Now we can define the non-pluripolar product of $\varphi_i \in {\rm PSH}(X,\theta^i)$, $i=1,...,p$. It has been shown in \cite{BEGZ10} that the sequence of currents
$$
\mathbf{1}_{\cap_i \{\varphi_i>V_{\theta^i}-k\}} (\theta^1+dd^c\max(\varphi_1,V_{\theta^1}-k))\wedge...\wedge (\theta^p+dd^c\max(\varphi_p,V_{\theta^p}-k))
$$
is non-decreasing in $k$ and converges weakly to the so called {\it non-pluripolar product}
$$
\theta^1_{\varphi_1}\wedge...\wedge\theta^p_{\varphi_p}.
$$
By \cite[Proposition 1.4]{BEGZ10}, we know that the non-pluripolar product is symmetric and multilinear.
When $p = n$, the resulting positive Borel measure, which does not charge pluripolar sets. Pluripolar sets are Borel measurable sets that are contained within some set $\{\psi = -\infty\}$, where $\psi \in {\rm QPSH}(X)$. We call a measure on $X$  that does not charge pluripolar sets a non-pluripolar measure.
For a $\theta$-psh function $\varphi$, the {\it non-pluripolar complex Monge--Ampère measure} of $\varphi$ is
$$
\theta_\varphi^n := \lim_{k\to\infty} \mathbf{1}_{\{\varphi>V_\theta-k\}} (\theta+dd^c\max(\varphi,V_\theta-k))^n.
$$

The following volume comparison is due to Witt Nyström. See also \cite{DDL18,Vu21} for more general version.
\begin{thm}$($\cite[Theorem 1.2]{WN19}$)$\label{thm WN}
    Let $\varphi$ and $\psi$ be two $\theta$-psh functions. If $\varphi$ is less singular than $\psi$, then
    $$
    \int_X \theta_\varphi^n \geq \int_X \theta_\psi^n.
    $$
\end{thm}

  \subsubsection{Envelope} In our study of relative pluripotential theory, the following envelope construction will be essential. 
  \begin{defi}\label{def 1}
  The rooftop envelope $P_\theta(\psi, \varphi)$ is defined by $P_\theta(\psi, \varphi) := P_\theta(\min \{\psi, \varphi\})$.
  Given $\psi,\varphi \in {\rm PSH}(X, \theta)$, envelopes with respect to singularity $P_\theta[\psi](\varphi)$ is defined by
  $$
  P_\theta[\psi](\varphi):= \left( \lim_{C\to +\infty} P_\theta(\psi+C,\varphi)\right)^*.
  $$
When $\varphi=V_\theta$, we simply write $P_\theta[\psi] := P_\theta[\psi](V_\theta)$.
\end{defi}
The envelopes with respect to singularity was introduced by Ross and Witt Nyström \cite{RW14}, building on ideas of Rashkovskii and Sigurdsson \cite{RS05} in the local setting.\par

A $\theta$-{\it model potential} $\phi$ is a $\theta$-psh function such that $\phi=P_\theta[\phi]$ and $\int_X \theta_\phi^n>0$. By definition, $V_\theta$ is clearly a $\theta$-model potential. The following theorem given by Darvas--Di Nezza--Lu \cite[Theorem 3.14]{DDL23} is frequently used in our proof.
\begin{thm}\label{thm 2.2}
    Assume that $\varphi \in {\rm PSH}(X,\theta)$ such that $\int_X \theta_\varphi^n>0$ and $\varphi \leq 0$. Then $\varphi \leq P_\theta[\varphi]$.
\end{thm}

\subsubsection{Relative full mass class}
\begin{defi}\label{def 2}
Given a potential $\phi\in {\rm PSH}(X, \theta)$ such that $\int_X \theta_\phi^n>0$. The relative full mass class is defined by 
$$
\mathcal{E}(X, \theta, \phi):= \left\{ u \in {\rm PSH}(X,\theta): u \preceq \phi,~\int_X \theta_u^n = \int_X \theta_\phi^n \right\}.
$$
The full mass class is defined by $\mathcal{E}(X,\theta):= \mathcal{E}(X,\theta,V_\theta)$.
\end{defi}

The following theorem comes from \cite[Theorem 3.15]{DDL23} and characterizes the relationship between $u \in \mathcal{E}(X,\theta,\phi)$ and $P_\theta[u]$.
\begin{thm}\label{thm 2.4}
    Assume that $\phi$ is a $\theta$-model potential. Then $u \in \mathcal{E}(X,\theta,\phi)$ iff $u \preceq \phi$ and $P_\theta[u] = \phi$.
\end{thm}
The following result, known as the uniqueness principle, states that the non-pluripolar Monge--Ampère measure determines the potential within a relative full mass class.
\begin{thm}$($\cite[Theorem 3.13]{DDL23}$)$\label{thm 2.5}
   Assume $\phi$ is a model potential and $v,u \in \mathcal{E}(X,\theta,\phi)$. Then
   $$
   \theta_u^n = \theta_v^n~ \Longleftrightarrow ~u\equiv v ~{\rm up~ to ~a ~constant.}
   $$
\end{thm}

 Let $x \in X$. Fixing a holomorphic chart $x \in U \subset X$, the Lelong number $\nu(\varphi,x)$ of $\varphi \in {\rm PSH}(X, \theta)$ at $x$ is defined as follows:
 $$
 \nu(\varphi,x) := \sup\{ a \geq 0 :~\varphi(z) \leq a \log |z-x| +O(1)\}.
 $$
The Lelong number $\nu(\varphi)$ of $\varphi \in {\rm PSH}(X, \theta)$ is defined by $\nu(\varphi) : = \sup_{x\in X} \nu(\varphi,x)$.

The following theorem is also from the work of Darvas–Di Nezza–Lu, which generalized a special case of \cite[Theorem 1.1]{BBEGZ} and solved an open problem (see \cite[Remark 10.3]{GZ17}).
\begin{thm}$($\cite[Theorem 1.1.(1)]{DDNL18}$)$\label{thm 2.7}
    For any $\varphi \in \mathcal{E}(X,\theta)$, one has
    $$
    \nu(\varphi,x) = \nu(V_\theta,x),~\forall x \in X.
    $$
\end{thm}
See also \cite[Lemma 5.1]{DDL23} for the relative version.

\subsubsection{Monge--Ampère Capacity} Let $\psi \in {\rm PSH}(X,\theta$). The relative  Monge–Ampère capacity of a Borel set $E \subset X$ is defined as
$$
{\rm Cap}_\psi(E) := \sup \left\{  \int_E \theta_u^n:~u \in {\rm PSH}(X,\theta),~\psi-1\leq u \leq \psi \right\}.
$$
The Monge–Ampère capacity is then given by ${\rm Cap}_\theta := {\rm Cap}_{V_\theta}$.
Note that the relative Monge--Ampère capacity is inner regular, one can see \cite[Lemma 4.3]{DDL23}.

\section{Proof of $L^1$-stability}\label{sec 3}
Assume that $\{ \theta \},\{ \theta^j \}$ for $j=1,2,...,$ are big cohomology classes. Let's consider the following complex Monge--Ampère equations:
\begin{equation}\label{eq 3.11}
\begin{split}
    (\theta^j_{\varphi_j})^n = \mu_j,&~{\rm and}~\sup_X \varphi_j  =0; j=1,2,...,\\
 \theta_\varphi^n = \mu,  &~{\rm and}~\sup_X \varphi  =0.
\end{split}
\end{equation}
We require that $\int_X \mu_j,\int_X \mu>0$ and set $\phi = P_\theta[\varphi]$, $\phi_j = P_{\theta^j}[\varphi_j]$. Then, we have
\begin{thm}\label{Proposition A'}
    Assume that $\theta^j \to_{\mathcal{C}^+} \theta$, and $\phi_j \to_{L^1(\omega^n)} \phi$. If $\| \mu_j - \mu  \| \to 0$, then
    $$
    \|  \varphi_j - \varphi \|_{L^1(\omega^n)} \to 0.
    $$
\end{thm}

Let $\theta$ be a smooth, closed real $(1,1)$-form on $X$ whose cohomology class is big. Assume $\epsilon_j \searrow 0$ and set $\alpha^j := \theta+ \epsilon_j \omega$. Then we have $\{ \alpha^j \}$ are big, $\alpha^j \searrow_{\mathcal{C}^0} \theta$, and ${\rm PSH}(X,\theta)\subset{\rm PSH}(X,\alpha^{j+1})\subset {\rm PSH}(X,\alpha^j)$.

The following lemma is a slight generalization of \cite[Lemma 5.16]{DDL23}, which yields the $L^1$-limit of subsolutions is also a subsolution, even if the $(1,1)$-forms $\alpha^j$ have slight variations.

\begin{lem}\label{lem 3.2}
    Let $u_j\in {\rm PSH}(X,\alpha^j)$ such that $(\alpha^j_{u_j})^n \geq f_j\mu$, where $ 0 \leq f_j \in L^1(X, \mu)$ and $\mu$ is a non-pluripolar measure on $X$. Assume that $f_j\to_{L^1(\mu)} f \in L^1(\mu)$, and $u_j \to_{L^1(\omega^n)} u \in {\rm PSH}(X,\theta)$. Then $\theta_u^n \geq f \mu$.
\end{lem}
\begin{proof}
    We fix $i \in \mathbb{Z}^+$. For all $j\geq i$, we have $u_j \in {\rm PSH}(X,\alpha^i)$. 
    Since we have $u_j \to_{L^1(\omega^n)} u \in {\rm PSH}(X,\alpha^i)$ and $f_j \to_{L^1(\mu)} f$ as $j \to \infty$, applying \cite[Lemma 5.16]{DDL23}, we get $(\alpha^i_{u})^n \geq f\mu$. 
    Note that the non-pluripolar product is multilinear, we obtain  
    $$
    (\alpha^i_{u})^n = \theta_u^n + \sum_{k=1}^n \epsilon_i^k \binom{n}{k}  \theta_u^{n-k}\wedge\omega^{k} =: \theta_u^n + m_i.
    $$
    It is easy to see that $0 \leq \int_X m_i = O(\epsilon_i)$, which implies $(\alpha^i_{u})^n \to \theta_u^n$ in the weak sense. Hence $\theta_u^n \geq f\mu$.
\end{proof}

\begin{prop}\label{prop 2.5}
    Let $\mu_j$, $\mu$ be non-pluripolar measures on $X$ such that $\sup_j \mu_j(X) <+\infty$. Then
    $$
    \| \mu_j - \mu \| \to 0
    $$
    if and only if there exists a non-pluripolar Radon measure $\nu$ on $X$ and $f,f_j \in L^1(\nu)$, such that $\mu_j := f_j\nu$, $\mu := f\nu$ and $f_j \to_{L^1(\nu)} f $.
\end{prop}
Proposition \ref{prop 2.5} follows from the argument in \cite[page 1034]{1}, which gives an equivalent expression of Theorem  \ref{Proposition A}. In the setting of Theorem \ref{Proposition A}, assume $ (\theta^j_{\varphi_j})^n = f_j \nu$ and $ \theta_\varphi^n = f\nu$, where $\nu$ is a non-pluripolar measure and $f_j,f\in L^1(\nu)$. If we have $f_j \to_{L^1(\nu)} f$, then 
$$
\| \varphi_j - \varphi \|_{L^1(\omega^n)} \to 0.
$$
\begin{proof}
 Assume $\|  \mu_j - \mu \| \to 0$. We define
 $$
  \nu:= \mu + \sum_{j\geq 1} 2^{-j} \mu_j.
 $$
  Since $\sup_j \mu_j(X) <+\infty $, we have $\nu$ is a well defined non-pluripolar Radon measure with respect to which $\mu_j,\mu$ are absolutely continuous. So, we have $\mu_j := f_j\nu$, $\mu := f\nu$, where $f_j,f\in L^1(\nu)$, and
 $$
   \|\mu_j - \mu\| = \| f_j - f \|_{L^1(\nu)}     
 $$
 by definition. Now we can obtain the conclusion directly. The converse is obvious by definition of the total variation of the signed measure.
\end{proof}

\begin{proof}[Proof of Theorem \ref{Proposition A'}]
 Since $\theta^j \to_{\mathcal{C}^+}\theta$, for all $\epsilon>0$, there exists $j_0>0$, such that
 $$
 \theta^j - \theta \leq \epsilon \omega
 $$
 for all $j>j_0$. Then we can find a suitable sequence $\epsilon_j \searrow 0 $ such that 
 $$
 \theta^j < \theta+ \epsilon_j \omega =: \alpha^j .
 $$
 Hence $\varphi_j \in {\rm PSH}(X,\alpha^j)$.\par
 Since $\{ \varphi_j \}_{j\geq k} \subset {\rm PSH}(X,\alpha^k)$ and $\sup_X \varphi_j = 0$. By weak compactness, we can assume-up to extracting-that 
 $$
 \varphi_j \to_{L^1(\omega^n)} \psi \in {\rm PSH}(X,\alpha^k), ~\forall  k >0.
 $$
 This yields $\psi \in {\rm PSH}(X,\theta)$ with $\sup_X \psi = 0$.\par
 \textit{Now we claim} $\psi = \varphi$.\\
 $1^\circ.$ By Theorem \ref{thm 2.2}, we have $\varphi_j \leq P_\theta[\varphi_j]$ and by Theorem \ref{thm 2.4}, we have $P_\theta[\varphi_j] =  \phi_j$. Since $\phi_j \to_{L^1(\omega^n)} \phi$, we obtain that $\psi \leq \phi$ almost everywhere with respect to $\omega^n$; hence everywhere because they are quasi-psh functions, see for example \cite[Corollary 1.38]{GZ17}.\\
 $2^\circ$. Set $\gamma_j:= \alpha^j - \theta^j >0$. Since the non-pluripolar product is multilinear, we have
\begin{equation}\label{eq 3}
(\alpha^j_{\varphi_j} )^n = \left(\theta^j_{\varphi_j}+ \gamma_j \right)^n = \mu_j + m_j \geq \mu_j,
\end{equation}
where $\mu_j = (\theta^j_{\varphi_j})^n $ and $m_j : = \sum_{k=1}^n {n \choose k} \gamma_j^k \wedge (\theta^j_{\varphi_j})^{n-k} \geq 0$. Applying Proposition \ref{prop 2.5}, there exists a non-pluripolar measure $\nu$, and $f_j,f\in L^1(\nu)$ such that
$$
\mu_j := f_j\nu,~\mu := f\nu ~{\rm and}~ f_j \to_{L^1(\nu)}f.
$$
Since we have $\varphi_j \to_{L^1(\omega^n)} \psi$, and (\ref{eq 3}), it follows from Lemma \ref{lem 3.2} that  
$$
\theta_\psi^n \geq f\nu = \mu.
$$
Note that we have $\psi \leq \phi$ by $1^\circ$. Using Theorem \ref{thm WN}, we obtain
$$
\int_X \theta_\psi^n \leq \int_X \theta_\phi^n = \int_X \theta_\varphi^n = \int_X \mu.
$$
Now, by comparing the total mass we thus obtain 
$$
\theta_\varphi^n=\mu=\theta_\psi^n,
$$
which implies that $\psi \in \mathcal{E}(X,\theta,\phi)$. Finally, by Theorem \ref{thm 2.5}, we have $\psi=\varphi$. Hence 
$$
\| \varphi_j - \varphi \|_{L^1(\omega^n)} = \| \varphi_j - \psi \|_{L^1(\omega^n)} \to 0.
$$
\end{proof}

By applying the same technique as in Theorem \ref{Proposition A'}, we derive the following two stability corollaries for (\ref{eq 3.11}):
\begin{cor}\label{Proposition B'}
    Assume that $\theta^j \to_{\mathcal{C}^+} \theta$, and $\phi_j \leq \phi$. If $\| \mu_j - \mu  \| \to 0$, then
    $$
    \|  \varphi_j - \varphi \|_{L^1(\omega^n)} \to 0.
    $$
\end{cor}
\begin{cor}\label{Proposition C'}
    Assume that $\theta^j \to_{\mathcal{C}^+} \theta$, and $\phi = V_\theta$. If $\| \mu_j - \mu  \| \to 0$, then
    $$
    \|  \varphi_j - \varphi \|_{L^1(\omega^n)} \to 0.
    $$
\end{cor}

\begin{proof}[Proof of Corollary \ref{Proposition B'}]
We continue to use the construction from the proof of Theorem \ref{Proposition A}. Through the weak compactness, we obtain that
$$
 \varphi_j \to_{L^1(\omega^n)} \psi \in {\rm PSH}(X,\alpha^k),~\forall k>0.
 $$
 This implies $\psi \in {\rm PSH}(X,\theta)$ and $\sup_X \psi =0$. By Theorem \ref{thm 2.2} and Theorem \ref{thm 2.4}, again, we have $\varphi_j \leq\phi_j$, which yields
 $$
 \varphi_j \leq \phi_j \leq \phi.
 $$
 Since $\varphi_j \to_{L^1(\omega^n)} \psi$, then $\psi \leq \phi$. The remaining part follows from the same steps as in the proof of Theorem \ref{Proposition A}.
\end{proof}

\begin{rem}
    Note that if we can compare the singularities of the $\theta$-model potential $\phi$ and the limiting function $\psi$ as above, such that $\phi \succeq \psi$, then it follows that $\int_X \theta_\phi^n \geq \int_X \theta_\psi^n$. By comparing the total mass, proving the stability of the solution is equivalent to prove the stability of the subsolution, namely, finding the variant of \cite[Lemma 5.16]{DDL23}. This is somewhat simpler than other methods. However, the limitation is that we can only achieve $L^1$-convergence for the solutions.
\end{rem}

Using the same method, we can directly obtain Corollary \ref{Proposition C'}. 
In conclusion, the crucial point is the stability of subsolutions, Lemma \ref{lem 3.2}. 
Thus, by extracting the proof of Theorem \ref{Proposition A'}, we establish the following generalization.
\begin{cor}
    Assume that $\{ \theta \},\{ \theta^j \}$, $j=1,2,...,$ are big cohomology classes on $X$ such that $\theta^j \to_{\mathcal{C}^+} \theta$. 
    Let $u_j\in {\rm PSH}(X,\theta^j)$ such that $(\theta^j_{u_j})^n \geq \mu_j$, where $\mu_j$ are non-pluripolar measures on $X$. If $\| \mu_j - \mu \|\to 0$ and $u_j \to_{L^1(\omega^n)} u \in {\rm PSH}(X,\theta)$, then $\theta_u^n \geq  \mu$.
\end{cor}
\begin{proof}
    Since $\theta^j \to_{\mathcal{C}^+} \theta$, we can find a suitable sequence $\epsilon_j \searrow0$ such that
    $$
   \theta^j < \theta+ \epsilon_j\omega =: \alpha^j.
    $$
    By the multilinearity of the non-pluripolar product, we obtain $(\alpha^j_{u_j}) \geq \mu_j$. 
    Then, applying Lemma \ref{lem 3.2} and Proposition \ref{prop 2.5}, we complete the proof.
\end{proof}

\section{Proof of $\mathcal{C}^{k,\alpha}$-stability}\label{sec 4}
\subsection{Di Nezza--Lu's estimate}\label{sec 4.1}
Let $\theta$ be a smooth, closed semi-positive real $(1,1)$-form on $X$ such that $\int_X \theta^n>0$. For the convenience, we briefly review some details from \cite{DNL14}. First, Di Nezza and Lu proved that

\begin{thm}$($\cite[Theorem 3]{DNL14}$)$\label{thm 4.1}
 Let $\psi^{\pm}$ be quasi-psh functions, such that $e^{- \psi^-} \in L^1(\omega^n)$
 and $\psi^{\pm} \in \mathcal{C}^\infty(X\backslash D)$, where $D$ is a divisor on $X$. Then there exists an effective s.n.c. divisor $E$ on $X$ such that $\{ \theta \} - c_1(E)$ is a Kähler class, and the solution $\varphi$ of the complex Monge--Ampère equation
 \begin{equation}\label{eq 4.0}
 \theta_\varphi^n = ce^{\psi^+ - \psi^-} \omega^n,~\varphi \in \mathcal{E}(X,\theta)
 \end{equation}
is smooth outside $D\cup E$, where $c$ is a normalized constant.
\end{thm}

\subsubsection{Demailly’s regularization}\label{sec 4.1.1}
In this subsection, we will borrow some basic notations from \cite{De94}. If $u \in {\rm PSH}(X,C\omega)$, set 
\begin{equation}\label{eq 4.22}
\rho_\epsilon(u)(p) := {\epsilon^{-2n}} \int_{\zeta \in T_{X,p}} u({\rm exph}_p(\zeta)) \chi\left( \epsilon^{-2}{|\zeta|^2} \right) d\lambda(\zeta),~\epsilon>0,
\end{equation}
then one has $\rho_\epsilon(u) \in \mathcal{C}^\infty(X)$. Here $0\leq \chi \in \mathcal{C}^\infty(\mathbb{R})$ is a cut-off function such that ${\rm supp}\chi \subset [-1,1]$, $\int_\mathbb{R} \chi =1$, and
\begin{align*}
{\rm exph}_p : T_{X,p}\to X
\end{align*}
is the formal holomorphic part of the Taylor expansion of the Exponential map defined by the metric $\omega$. We rewrite $\rho_i(u)$ as $ \rho_{1/i}(u)$. Then we have $\rho_i(u) \searrow u$ as $i\to \infty$, and 
$$
dd^c\rho_i(u) \geq -(C+ \lambda_i)\omega,
$$
where $\lambda_i \searrow  \nu(u)$ as $i \to \infty$.

\subsubsection{Tsuji's trick}\label{sec 4.1.2} The so-called Tsuji’s trick \cite{Tsu} played a crucial role in the regularization of the degenerate complex Monge--Ampère equation (\ref{eq 4.0}).
That is, by adding a term $i^{-1}\omega$ to the $(1,1)$-form $\theta$, one can reduce the equation (\ref{eq 4.0}) to the well-known non-degenerate case.
We normalized $f := e^{\psi^+ - \psi^-}$ such that $\int_X f\omega^n = \int_X \theta^n$. Consider the complex Monge--Ampère equations
\begin{equation}\label{eq 4.1}
\left(\theta+ \frac{1}{i} \omega + dd^c \varphi_i' \right)^n = c_i' e^{\rho_i(\log f)} \omega^n,~\sup_X \varphi_i' =0,
\end{equation}
where $c_i'$ are normalized constants.
It follows from \cite{Yau} that there exists a unique solution $\varphi_i' \in {\rm PSH}(X,\theta+i^{-1}\omega) \cap \mathcal{C}^\infty(X)$ to (\ref{eq 4.1}) for each $i \in \mathbb{Z}^+$. 
 Since $e^{\rho_i(\log f)}\to f$ pointwise and $e^{\rho_i(\log f)} \leq\rho_i(f)$, by Jensen’s inequality, the Dominated Convergence Theorem yields $ e^{\rho_i(\log f)} \to_{L^1(\omega^n)} f $. So
$$
\int_X \left(\theta+ \frac{1}{i} \omega + dd^c \varphi_i' \right)^n = \int_X \left(\theta+ \frac{1}{i} \omega \right)^n \to \int_X \theta^n = \int_X \theta_\varphi^n= \int_X f\omega^n
$$
   gives that $c_i' \to 1$. Then, by Theorem \ref{Proposition A} we have $\varphi_i' \to_{L^1(\omega^n)}\varphi$.\par

  Hence if $\{\varphi_i'\}$ is pre-compact for the $\mathcal{C}^{k,\alpha}_{loc}(X\backslash(D\cup E))$-topology, $\forall k \in \mathbb{Z}^+$, then one has 
  $$
  \varphi_i' \to_{\mathcal{C}^{k,\alpha}_{loc}(X\backslash(D\cup E))} \varphi,~\forall k>0.
  $$
  Due to the Evans--Krylov method and Schauder estimates (see \cite[Chapter 14.3]{GZ17} for example), we only need to deal with the $\mathcal{C}^0$ estimate and the Laplace estimate of $\varphi_i'$ on $X \backslash (D\cup E)$.

  \subsubsection{Uniform Skoda Integrability Theorem} It follows from Skoda’s theorem that $\int_X e^{- C u}\omega^n$ is finite for all $C < 2/\nu(u)$, where $u \in {\rm QPSH}(X)$.
  There is a uniform integrability result, directly from Zeriahi \cite{Ze01}:
  \begin{lem}\label{lem 4.2}
      Let $\mathcal{U}\subset{\rm PSH}(X,A \omega)$ be a compact family for the $L^1(\omega^n)$-topology. Set $\sup_{u\in \mathcal{U}} \nu(u ) = a \geq 0$. Then there exists $C_2:= C_2 \left( C_1, \mathcal{U},A,\omega \right)$ such that
      $$
      \sup_{ u \in \mathcal{U}} \int_X e^{-C_1 u } \omega^n \leq C_2
      $$
      for all $0< C_1 < 2/a$.
  \end{lem}

\subsubsection{$\mathcal{C}^0$ estimate and Laplace estimate}\label{sec 4.1.4} Since $\theta\geq 0$ and $\int_X \theta^n >0$, by Kodaira's Lemma, there exists an effective  s.n.c. $\mathbb{R}$-divisor $E = \sum_i^N a_iE_i$ such that $ \alpha :=\{ \theta \} -c_1(E) $ is a Kähler class. Let  $\omega_0 = \theta - c_1(E,h) \in \alpha$ be a Kähler form. Set 
$$
\Phi: = \sum_i^N a_i \log |\sigma_i |_{h_i},
$$
where $\sigma_i \in \mathcal{O}(E_i)$ such that $\sigma_i$ vanishes on $E_i$.
Then we have
\begin{equation}\label{eq 4.2}
\theta +dd^c \Phi = \omega_0 + [E].
\end{equation}
by the Poincaré--Lelong equation. By rescaling $\omega$ we may assume that $\omega_0 \geq \omega$.\\

 Now, we introduce results from \cite[Theorem 5.1, 5.2]{DNL14}.
\begin{thm}$(\mathcal{C}^0 ~{\rm estimate})$\label{thm 4.3'}
    Assume $f= e^{\psi^+ - \psi^-} \in {L}^1(\omega^n)$ such that $\psi^\pm \in {\rm QPSH}(X) \cap L^\infty_{loc}(X\backslash D)$. Let $C>0$ such that $\psi^\pm \in {\rm PSH}(X,C\omega)$, $\sup_X \psi^+ \leq C$, and $\theta \leq C\omega$. Let $\varphi$ be the unique normalized solution to
$$
  \theta_\varphi^n = c f \omega^n,~ \varphi \in \mathcal{E}(X,\theta)
$$
Then, for any $a > 0$ satisfies $a\psi^- \in  {\rm PSH}(X,{\omega_0}/{2})$, there exists a constant $A = A\left( C,\omega,\int_X e^{-2\varphi/a}\omega^n \right) > 0$ such that
$$
\varphi \geq a\psi^- + \Phi -A.
$$
\end{thm}

\begin{thm}$($Laplace estimate$)$\label{thm 4.3}
Let $f= e^{\psi^+ - \psi^-}$, where $\psi^\pm \in \mathcal{C}^\infty(X)$. Fix $t\in(0,1)$. Let $\varphi \in \mathcal{C}^\infty(X) \cap {\rm PSH}(X,\theta+ t \omega)$  be the unique normalized solution to
$$
(\theta+ t\omega + dd^c \varphi)^n = c_t e^{\psi^+ - \psi^-}\omega^n,~\sup_X \varphi = 0.
$$
Assume given a constant $C > 0$ such that
$$
dd^c \psi^\pm \geq -C\omega,~\sup_X \psi^+ \leq C,~{\rm and}~\theta \leq C\omega.
$$
Assume also that the holomorphic bisectional curvature of $\omega$ is bounded from below by $-C$. Then there exists a constant $A:= A\left( C, \omega ,\int_X e^{-2(3C+1)\varphi }\omega^n \right) > 0$ such that
$$
\Delta_\omega \varphi  = {\rm tr}_\omega dd^c \varphi \leq   A\cdot e^{-\psi^- - (3C+1)\Phi}.
$$
\end{thm}

Go back to (\ref{eq 4.1}), note that we have $\mathcal{U}:= \{ \varphi_i' \}\cup \{ \varphi \}$ is compact for the $L^1(\omega^n)$-topology. Then, by Lemma \ref{lem 4.2} we have
$$
\int_X e^{-C_1 \varphi_i' }\omega^n \leq C_2(C_1, \mathcal{U},C,\omega),~\forall~C_1>0.
$$
 We thus obtain $\varphi_i',~\Delta_\omega \varphi_i'$ are locally uniformly bounded outside $D\cup E$.

\subsection{$\mathcal{C}^{k,\alpha}$-stability}\label{sec 4.2}
Let $\theta^j$ be smooth, closed semi-positive real $(1,1)$-forms on $X$ such that $\int_X (\theta^j)^n>0$ for all $j$. Assume that $\theta^j \to_{\mathcal{C}^0} \theta$. Then we may assume that
$$
- \frac{1}{2} \omega_0 \leq \theta^j -\theta \leq \frac{1}{2} \omega_0,
$$
after possibly discarding a finite number of terms. We then obtain $\{ \theta^j\} - c_1(E)$ are Kähler classes, where $E$ and $\omega_0$ are defined in section \ref{sec 4.1.4}. We get
$$
\theta^j + dd^c \Phi = \omega^j + [E],
$$
where $\omega^j := \theta^j - c_1(E,h) \geq \frac{1}{2}\omega_0$. By rescaling $\omega$, we can assume that $\omega^j, \omega_0 \geq \omega$ for all $j$, and we fix $\omega$ afterwards.
 Applying Theorem \ref{thm 4.1} and the arguments in section \ref{sec 4.1}, that the solutions $\varphi_j \in \mathcal{E}(X,\theta^j)$ of the complex Monge--Ampère equations
\begin{equation}\label{eq 4.3}
    (\theta^j_{\varphi_j})^n = c_j e^{\psi^+ - \psi^-} \omega^n,~\sup_X \varphi_j = 0, 
\end{equation}
are smooth outside $D\cup E$. 
\begin{thm}\label{thm D'}
    Let $U\Subset X \backslash (D\cup E)$ be an arbitrary holomorphic coordinate chart, such that $\theta^j = dd^c g_j$, $\theta = dd^c g$ on $U$. Assume that $g_j \to_{\mathcal{C}^{\infty}(U)} g$, then   
    $$
    {\varphi_j} \to_{\mathcal{C}^{\infty}_{loc}(X \backslash(D \cup E))} \varphi.
    $$
\end{thm} 


\begin{rem}
    If one can find smooth, closed real $(1,1)$-forms $\alpha^i,i=1,...,N$ such that 
    $$
    \theta^j = \sum_{i} t_{ij} \alpha^i,~\theta = \sum_i t_i \alpha^i, ~{\rm and} ~\lim_j t_{ij} = t_i,
    $$
    then the assumptions of Theorem \ref{thm D'} are trivially satisfied.
\end{rem}

Set $f = e^{\psi^+ - \psi^-}$. Consider the  complex Monge--Ampère equations
\begin{equation}\label{eq 4.4}
\left(\theta^j + \frac{1}{i}\omega + dd^c\varphi_{ij} \right)^n = c_{ij} e^{\rho_i(\log f)}\omega^n,~\sup_X \varphi_{ij}=0,
\end{equation}
where $c_{ij}$ are normalized constants. 
By \cite{Yau}, there exists a unique solution $\varphi_{ij} \in {\rm PSH}(X,\theta^j+i^{-1}\omega) \cap \mathcal{C}^\infty(X)$ to (\ref{eq 4.4}) for each $i,j$.
By arguments in section \ref{sec 4.1.1}, one may assume that
$$
dd^c \rho_i(\psi^\pm) \geq -C\omega, ~ \sup_X \rho_i(\psi^+) \leq C.
$$
Assume also that the holomorphic bisectional curvature of $\omega$ is bounded from below by $-C$, and that $\theta^j,\theta \leq C\omega$. 
Set $\mathcal{U}_j := \{ \varphi_{ij} \}_i \cup \{ \varphi_j \}$. Since 
\begin{align*}
    \int_X \left(\theta^j + \frac{1}{i}\omega + dd^c\varphi_{ij} \right)^n = \int_X &\left(\theta^j + \frac{1}{i}\omega  \right)^n\\
    &\to \int_X (\theta^j)^n = \int_X (\theta^j_{\varphi_j})^n ,
\end{align*}
and $e^{\rho_i(\log f)} \to_{L^1(\omega^n) } f $, it follows that  $\lim_ic_{ij}\to c_j$. Applying Theorem \ref{Proposition A}, $\mathcal{U}_j$ is $L^1(\omega^n)$ compact. By arguments in section \ref{sec 4.1.2}, we see that $\varphi_{ij} \to_{\mathcal{C}^\infty_{loc}(X\backslash D\cup E))} \varphi_i$ as $j\to \infty$. 
Then, by using Lemma \ref{lem 4.2} and Theorem \ref{thm 4.3}, we obtain
\begin{equation}\label{eq 4.3'}
\Delta_\omega \varphi_j \leq A_j e^{-\psi^- - (3C+1)\Phi},
\end{equation}
where $A_j := A\left(C,\omega \right )\cdot  \sup_{u\in\mathcal{U}_j} \int_X e^{-2(3C+1)u}\omega^n $. 
Set $\mathcal{U}' := \left( \cup_j \mathcal{U}_j \right) \bigcup \mathcal{U}$, where $\mathcal{U} = \{ \varphi_i' \} \cup \{\varphi\}$ (see section \ref{sec 4.1.2}), we have the following result:
\begin{prop}\label{cor 4.4}
    The family $\mathcal{U}' \subset {\rm PSH}(X, (C+1) \omega)$ is compact for the $L^1(\omega^n)$-topology, and $\nu( u )= 0 $, for all $u \in \mathcal{U}'$.
\end{prop}
\begin{proof} For $\forall u \in \mathcal{U}'$, we have
$$
u \in \mathcal{E}\left(X,\theta^j + \frac{1}{i} \omega\right)~{\rm or}~u \in \mathcal{E}\left(X,\theta+ \frac{1}{i} \omega \right)
$$
for some $i \in \mathbb{Z}^+ \cup \{\infty\},j \in \mathbb{Z}^+$. Then we obtain $\nu( u ) = 0$ by Theorem \ref{thm 2.7}. Now let's deal with compactness. Set $\{ u_k \}\subset \mathcal{U}'$:\\
    $1^\circ.$ Suppose that the intersection $\{ u_k \} \cap \mathcal{ U}$ is infinite. Then we may assume that $\{ u_k \} \subset \mathcal{U}$. Since $\mathcal{U}$ is compact, we have $u_k \to_{L^1(\omega^n)} \varphi$.\\
    $2^\circ.$ We may assume that $\{ u_k \} \subset \{  \varphi_{ij}\}_{i,j}$. Set $u_k = \varphi_{i_k,j_k}$, we can find a subsequence of $\{ \varphi_{i_k,j_k} \}$, rewritten as $\{ \varphi_{i_k,j_k} \}$ such that $i_k\to\infty,j_k \to \infty$ as $k\to\infty$. Consider the complex Monge--Ampère equations
    $$
    \left( \theta^{j_k} + \frac{1}{i_k}\omega + dd^c \varphi_{i_k,j_k} \right)^n = c_{i_k,j_k} e^{\rho_{i_k}(\log f)}\omega^n,
    $$
    and
    $$
    \left( \theta + dd^c \varphi \right)^n =  f\omega^n.
    $$
    Since $\theta^{j_k} + \frac{1}{i_k}\omega \to_{\mathcal{C}^0} \theta$, we have
    \begin{align*}
     \int_X \left( \theta^{j_k} + \frac{1}{i_k}\omega + dd^c \varphi_{i_k,j_k} \right)^n =& \int_X  \left( \theta^{j_k} + \frac{1}{i_k}\omega\right)^n \\
     &\to \int_X \theta^n =
     \int_X \left( \theta + dd^c \varphi \right)^n.
    \end{align*}
    We also have $ e^{\rho_{i_k}(\log f)} \to_{L^1(\omega^n)}  f$ as $k \to \infty$, which yields $ c_{i_k,j_k} \to 1$ as $k\to\infty$. By Theorem \ref{Proposition A} again, we get
    $$
    \varphi_{i_k,j_k} \to_{L^1(\omega^n) } \varphi~{\rm as}~k \to \infty.
    $$
    $3^\circ.$ There are three other situations: fix $i_0$, suppose that the intersection $\{ u_k \} \cap \{ \varphi_{i_0j}\}_j$ is infinite; or fix $j_0$, suppose that the intersection $\{ u_k \} \cap \mathcal{U}_{j_0}$ is infinite; or suppose that the intersection $\{  u_k \}_k \cap \{ \varphi_j \}_j$ is infinite.
    Repeating the above argument for each case, we complete the proof.
\end{proof}

\begin{proof}[Proof of Theorem \ref{thm D'}]
     Since $\omega^j \geq \frac{1}{2}\omega_0$, we can find a uniform constant $a$ such that $a\psi^- \in \mathrm{PSH}(X, \omega_0/4) \subset \mathrm{PSH}(X, \omega^j/2)$. By Proposition \ref{cor 4.4} and Lemma \ref{lem 4.2} we have $\sup_{u \in \mathcal{U}'} \int_X e^{-2u/a} \omega^n<+\infty$. It follows from Theorem \ref{thm 4.3} and \ref{thm 4.3'}, there exists $A>0$ such that
     $$
     \varphi_j \geq a \psi^- + \Phi -A,~\Delta_\omega \varphi_j \leq Ae^{-\psi^- -(3C+1)\Phi}.
     $$
    Let $U\Subset X\backslash (D\cup E)$ such that $\theta^j =dd^c g_j$, $\theta = dd^c g$ on $U$.
    We have $\varphi_j,\varphi\in \mathcal{C}^\infty(U)$ due to Theorem \ref{thm 4.1}. By the above argument, $g_j + \varphi_j$ and $\Delta_\omega (g_j + \varphi_j)$ are uniformly bounded on $U$.
     Let $U' \Subset U$. Applying Evans--Krylov method and Schauder estimates, we get
    $$
    \left\| g_j + \varphi_j \right\|_{\mathcal{C}^{k ,\alpha}(U')} \leq C_{U',k,\alpha},
    $$
    Since we have $g_j \to_{\mathcal{C}^{k ,1}(U)} g$, we thus obtain $\{ \varphi_j \}$ is pre-compact for the $\mathcal{C}^{k,\alpha}(U')$-topology, and any limit point $\psi$ of $\{ \varphi_j\}$ satisfies $\psi = \varphi$ almost everywhere on $U'$ by Theorem \ref{Proposition A}, hence everywhere. This implies that
    $$
    \varphi_j \to_{\mathcal{C}^{k,\alpha}(U')} \varphi,~\forall k \in \mathbb{Z}^+.
    $$
    Since $U,U'$ are arbitrary, the proof is done.
\end{proof}

\subsection{Further discussion}\label{sec 4.3}

Actually, we can provide a slightly general version of Theorem \ref{thm D'}. 
We still use the definitions and notations from the proof of Theorem \ref{thm D'} in this subsection.

\subsubsection{Local $\mathcal{C}^{k,\alpha}$-stability}\label{sec 4.3.1}
Let $ \psi_j^\pm  \in {\rm PSH}(X,C\omega)\cap \mathcal{C}^\infty(X\backslash D)$ such that $e^{-\psi_j^-} \in L^1(\omega^n)$.
Consider the complex Monge--Ampère equations
\begin{equation}\label{eq 4.90}
(\theta^j_{\varphi_j})^n = c_j e^{\psi^+_j - \psi^-_j} \omega^n,~\varphi_j \in \mathcal{E}(X,\theta^j),~\sup_X \varphi_j = 0.
\end{equation}
 Since $ \theta^j  - c_1(E,h)$ are Kähler forms (see section \ref{sec 4.1.4}), we have $\varphi_j \in \mathcal{C}^\infty(X \backslash (D \cup E))$ by Theorem \ref{thm 4.1}.\par

 Set $f_j : = e^{\psi^+_j - \psi^-_j}, f:= e^{\psi^+ - \psi^-}$, and assume that $f_j \to_{L^1(\omega^n)} f$. 
Using Tsuji's trick, we consider instead the following complex Monge--Ampère equations
\begin{equation}\label{eq 4.99}
\left(\theta^j + \frac{1}{i}\omega + dd^c\varphi_{ij} \right)^n = c_{ij} e^{\rho_i(\log f_j)}\omega^n,~\sup_X\varphi_{ij} =0,
\end{equation}
where $c_{ij}$ are normalized constants.
For each $i,j$, the equation (\ref{eq 4.99}) has a unique solution $\varphi_{ij}\in {\rm PSH}(X,\theta^j+i^{-1}\omega) \cap \mathcal{C}^\infty(X)$.
Since $\psi_j^\pm \in {\rm PSH}(X,C\omega)$, we obtain $\sup \nu(\psi_j^\pm) <  + \infty$. Then, we may assume that
$$
dd^c \rho_i(\psi^\pm_{j_0}) \geq -C\omega,
$$
after possibly discarding finitely many of $\left\{ \rho_i (\psi^\pm_{j_0}) \right\}_i$ for each $j_0$, and adjust $C$. 
If we further require that
\begin{enumerate}
  \item[(R1)] $\psi^+,\psi_j^+ \leq C$;
  \item[(R2)] $\psi^-,\psi_j^- \geq h$, where $h: X \to \mathbb{R} \cup \{ -\infty\}$ such that $e^{-h} \in L^1(\omega^n)$.
\end{enumerate}

Since $f_j \to_{L^1(\omega^n)} f$, we can find a subsequence $f_{j_k}$ such that $f_{j_k} \to f$ a.e. as $k\to \infty$. By $f_j,f \in \mathcal{C}(X\backslash D)$, we deduce that $\log f_{j_k} \to \log f $ a.e. as $k \to \infty$. Combined with $\rho_i (\log f_j) \to \log f_j$ point-wise as $i \to\infty$ for $\forall j$, we can choose a subsequence $\{ i_k,j_k \}$ of $\{i,j\}$ such that $\rho_{i_k}(\log f_{j_k}) \to \log f$ a.e. as $k\to \infty$. It follows from dominant convergence theorem that we have 
$$
e^{\rho_{i_k}(\log f_{j_k})} \to_{L^1(\omega^n)} f~{\rm as}~k \to \infty.
$$
 So, we can derive Proposition \ref{cor 4.4} by the same arguments. Assume also that the holomorphic bisectional curvature of $\omega$ is bounded from below by $-C$, and that $\theta^j,\theta \leq C\omega$.
 Applying Theorem \ref{thm 4.3'} and Theorem \ref{thm 4.3} we get
\begin{equation}\label{eq 4.133}
\varphi_j \geq a \psi_j^- + \Phi - A~{\rm and}~\Delta_\omega \varphi_j \leq Ae^{-\psi^-_j - (3C+1)\Phi}.
\end{equation}
  By the same arguments in the proof of Theorem \ref{thm D'}, we have
\begin{prop}\label{thm 4.9}
    Let $U \Subset X \backslash (D \cup E)$ be an holomorphic coordinate chart such that $\theta^j = dd^c g_j$, $\theta = dd^c g$ on $U$. 
    If $g_j \to_{\mathcal{C}^{k+2,\beta_1}(U)} g$ and  $ \psi_j^+ - \psi_j^- \to_{\mathcal{C}^{k,\beta_2}(U)} \psi^+ - \psi^- $, for some $k \in \mathbb{Z}^+$ and $\beta_1,\beta_2 \in (0,1)$, and the functions $\psi_j^\pm$ satisfy requirements (R1,R2), then there exists $0<\alpha<1$ such that
    $$
    \varphi_j \to_{\mathcal{C}^{k+2,\alpha}_{loc}(U)} \varphi.
    $$
\end{prop}


\begin{rem}
        In subsection \ref{sec 4.3.1}, we considered the equations (\ref{eq 4.90}) with $\psi^\pm_j \in {\rm PSH}(X, C\omega)$ and $e^{-\psi^-_j} \in L^1(\omega^n)$. By Guan–Zhou's strong openness theorem \cite{GZh15}, we actually have $e^{-\psi^-_j} \in L^{p_j}(\omega^n)$ for some $p_j > 1$. Hence, we have $\varphi_j \in L^\infty(X)$ by \cite[Theorem 1.3]{GZ07}. 
        However, we cannot guarantee that $\inf_j p_j > 1$, and meanwhile $\theta^j$ varies, 
        therefore, it is difficult to establish the uniform estimate of $\| \varphi_j \|_{L^\infty(X)}$.
        Nevertheless, this condition is unnecessary for Proposition \ref{thm 4.9}, due to an estimate given by \cite{DNL14}.
    \end{rem}

\subsubsection{The Kähler case}
In \cite[Section 3]{DNL14}, the authors also studied the complex Monge--Ampère equation on quasi-projective varieties (\ref{eq 4.0}) for $\theta>0$.
 They proved that the solution is smooth outside $D$, \cite[Theorem 1]{DNL14}.\par
 
 Let $ \psi_j^\pm  \in {\rm PSH}(X,C\omega)\cap \mathcal{C}^\infty(X\backslash D)$ such that $e^{-\psi_j^-} \in L^1(\omega^n)$. Consider the complex Monge--Ampère equations
 $$
 \left(\theta^j_{\varphi_j}\right)^n = c_j e^{\psi^+_j - \psi^-_j}\omega^n,~\varphi_j \in \mathcal{E}(X,\theta^j);
 $$
 where $\theta^j,j=1,2,...$ are Kähler forms, and $c_j$ are normalized constants. By \cite[Theorem 1]{DNL14}, we have $\varphi_j \in \mathcal{C}^\infty(X \backslash D)$. As an application of Proposition \ref{thm 4.9}, we have the following stability result under the Kähler setting.
\begin{cor}
    Assume $\theta^j$ satisfies $\{\theta^j\} \to \{ \theta \}$, and $ \psi_j^\pm $ satisfy requirements (R1,R2) and $\psi^+_j - \psi^-_j \to_{\mathcal{C}^\infty_{loc}(X\backslash D)} \psi^+ - \psi^-$, $e^{\psi^+_j - \psi^-_j} \to_{L^1(\omega^n)} e^{\psi^+ - \psi^-}$. Then we have
    $$
    \theta^j_{\varphi_j} \to_{\mathcal{C}^\infty_{loc}(X\backslash D)} \theta_\varphi.
    $$
    
      {\rm Here~$\{ \theta^j \} \to \{ \theta \}$~denotes~the~convergence~in~the~vector~space~$H^{1,1}(X,\mathbb{R})$}.
\end{cor}
\begin{proof}
We can find smooth, closed real $(1,1)$-forms $\beta_i,i=1,...,N$ such that $\theta = \sum_i  \beta_i$ and $\gamma^j = \sum_i t_{ij}\beta_i \in [\theta^j] $, where $\lim_j t_{ij} = 1$. 
We may assume $\gamma^j $ are Kähler forms.
Then, consider the solutions to the complex Monge--Ampère equations
$$
\left( \gamma^j_{\psi_j} \right)^n = c_j e^{\psi^+_j - \psi^-_j} \omega^n, ~ \psi \in \mathcal{E}(X,\gamma^j),
$$
we know that $\psi_j \in \mathcal{C}^\infty(X\backslash D) $ by \cite[Theorem 1]{DNL14}. It follows from Proposition \ref{thm 4.9} that $\psi_j \to_{\mathcal{C}^\infty_{loc}(U)} \varphi$ for all $U \Subset X \backslash D$. We then obtain
$$
\gamma^j_{\psi_j} \to_{\mathcal{C}^\infty_{loc}(X\backslash D)} \theta_\varphi.
$$
By Theorem \ref{thm 2.5}, we obtain 
$\theta^j_{\varphi_j}  =  \gamma^j_{\psi_j} $.
\end{proof}

\subsection{Application}\label{sec 4.4}

In this section, we briefly outline the application of the main theorem from Section \ref{sec 4} to the Calabi--Yau varieties. Following the construction in \cite[section 7.2]{EGZ09}, we define a Calabi–Yau variety $V$ to be a projective variety with canonical singularities and $K_V \sim_\mathbb{Q} 0$. According to \cite[Theorem 7.5]{EGZ09}, for an ample $\mathbb{R}$-divisor 
$D$ on $V$, there exists a unique singular Ricci-flat Kähler metric $\omega \in c_1(D)$ with bounded potential, and smooth on $V_{reg}$ satisfies $\omega^n = c \Omega \wedge \overline{\Omega}$.

When $D \in N^1(V)_{\mathbb{R}}$ is a nef and big divisor, let $D_j$ be a sequence of ample divisors converging to $D$, and denote by $\omega_j \in c_1(D_j)$ the corresponding singular Ricci-flat Kähler metrics. An interesting problem is to study the asymptotic behavior of $\omega_j$.

In \cite{To09}, Tosatti showed that when $V$ is a Calabi--Yau projective manifold, the sequence $\omega_j$ mentioned above converges smoothly outside a subvariety to a singular Ricci-flat Kähler metric. Inspired by Tosatti's work, we will consider the case where $V$ is $\mathbb{Q}$-factorial.

In this case, by \cite[Theorem 5.7]{Ka}, we know that ample cone of $N^1(V)_{\mathbb{R}}$ is polyhedral near $D$. Then, by Tosatti's arguments in \cite[Page 764]{To09}, we can find finitly many nef and big divisors $D_i'$ such that both $D$ and $D_j$ can be written as linear combinations of the $D_i'$ with the positive coefficients, and the coefficients of $D_j$ converging to those of $D$. Moreover, by \cite[Theorem 7.1]{HM}, the nef and big divisor on $V$ is semiample.

Let $\pi:X \to V$ be a log resolution, by \cite[Lemma 3.2]{BBEGZ}, the problem reduces to the stability of equation (\ref{eq 4.3}). Here $\theta_\varphi \in c_1(\pi^*D)$ and $\theta^j_{\varphi_j} = \pi^* \omega_j$ are the pull-backs of singular Ricci-flat Kähler metrics on $V$, and $\theta^j = \sum_i t_{ij} \theta_i'$, $\theta = \sum_i t_i \theta_i'$, where $0\leq \theta_i'\in c_1(\pi^* D_i')$ and $\lim_j t_{ij} = t_i$, $0 \leq t_i, t_{ij} $. 

Therefore, we can conclude that there exists a singular Ricci-flat Kähler metric $\omega_0 \in c_1(D)$($\theta_\varphi = \pi^* \omega_0$) with bounded potential, smooth outside $V_{sing}$ and a subvariety $Y$ such that $\omega_j \to \omega_0$ in the sense of currents and $\omega_j \to_{\mathcal{C}^\infty_{loc} (V_{reg} \backslash Y)} \omega_0$.

\begin{rem}
    In \cite{CT15,DGZ23}, they also considered generalizations of Tosatti’s result \cite{To09}. Moreover, when studying families of complex Monge–Ampère equations, their approach does not require the positivity of the $(1,1)$-forms (although their cohomology classes are assumed to be positive). In particular, \cite{DGZ23} also treats the case of singular Calabi--Yau varieties.
\end{rem}

\end{document}